\documentclass[11pt]{article}
\usepackage{latexsym,color,amsmath,amsthm,amssymb,amscd,amsfonts,graphicx,epstopdf,epsfig}
\usepackage{colortbl}

 \definecolor{forestgreen(traditional)}{rgb}{0.0, 0.27, 0.13}

\newcommand{\eps}{\varepsilon}
\newtheorem{theorem}{Theorem}
\newtheorem{lemma}{Lemma}
\newtheorem{corollary}[lemma]{Corollary}

\newtheorem{proposition}[lemma]{Proposition}
\newtheorem{definition}{Definition}

\def\R{{\mathbb R}}

\def\eps{{\varepsilon}}
\def\diam{{rm diam}}

\def\diam{\rm diam}
\def\conv{{\rm conv}}
\def\vol{{\rm vol}}

\begin{document}
\title {\Large\bf The isotropy constant and boundary properties of convex \vspace{2mm}\\
bodies
}
\author{Mathieu Meyer and Shlomo Reisner}

\date{}

\maketitle

\date{}

 \maketitle
 \begin{abstract} \noindent Let ${\cal K}^n$ be the set of all convex bodies in $\R^n$ endowed with the Hausdorff distance. We prove that if  $K\in {\cal K}^n$
has positive generalized Gauss curvature at some point of its boundary, then $K$ is not a local maximizer for the isotropy constant $L_K$.

\end{abstract}

\section{Introduction and statement of the main result.}
\vskip 2mm\noindent
Let $K$ be a
convex body in $\R^n$ endowed with its canonical scalar product and Euclidean norm denoted by $|\cdot|$. It is well known (as a standard reference to the subject we refer to
 \cite{BGVV}; another, earlier, comprehensive reference is \cite{MP}) that
there exists a unique (up to orthogonal
transformations) affine, volume preserving,  mapping $A:\R^n\to
\R^n$ such that
for some constant $M_K>0$, depending on $K$, one has for every $y\in \R^n$
$$\int_{AK} \langle x,y\rangle dx  =0\hbox{ and }
\int_{AK}\langle x,y\rangle^2 dx = M_K^2|y|^2.$$

\noindent We say that $K$ is {\it in isotropic position\/} (or that {\it $K$ is isotropic\/})
if $A$ is the identity on $\R^n$. The {\em isotropy constant\/}
$L_K$ of $K$ is defined by
$$L_K=\frac{M_K}{|K|^{\frac{n+2}{2n}}}\,.$$
where $|B |$ denotes the volume of a Borel subset $B$ of $\R^n$.
Note that it is customary to assume, as part of the definition of isotropic position, that $|AK|=1$; for the sake of convenience in our proofs, we prefer not to include this assumption in the definition.

\noindent
The famous {\em Slicing Problem} asks whether there exists a universal constant $C>0$ such that, for any $n$,  any convex body $K$ in $\R^n$ has a hyperplane section $K\cap H$ such that
$$\vol_{n-1}(K\cap H)\ge C\,\vol_n(K)^{\frac{n-1}{n}}.$$
This problem is equivalent to the existence of an upper bound $D>0$ for $L_K$, independent
of the dimension.
J. Bourgain proved in \cite{B} that $L_K\le Cn^{1/4} \log(n)$, this bound  was improved by B. Klartag in \cite{K}
to  $L_K\le Cn^{1/4} $, where $C$ is an absolute constant. Note that the minimum of $L_K$ is obtained only for ellipsoids (for an interesting discussion of stability in that inequality, see \cite{AB}).

\noindent Since the exact upper bound for $L_K$ is still an open problem, it is interesting to investigate
what are the properties of the maximizers for this quantity (a compactness argument shows that, for a fixed $n$,
maximizers for $L_K$ exist among convex bodies in $\R^n$).
We say that a convex body $K$ in $\R^n$ is a {\em local maximizer\/} (resp. {\em local minimizer}) for $L_K$ if for some
$\varepsilon>0$ one has
$L_{K'}\le L_K$  (resp. $L_{K'}\ge L_K$) for all convex bodies $K'$ in $\R^n$ such that $d(K',K)<\varepsilon$ ($d$ may
denote here the Hausdorff or the Banach--Mazur distance).
L. Rademacher proved in \cite{R} that if a simplicial polytope is a maximizer for $L_K$, then it  must be a simplex.
Campi, Colesanti and Gronchi
showed in \cite{CCG}, using shadow movements, that if $K$ has an open subset of its boundary which is
$C^2$ with positive Gauss curvature,
then $K$ can not be a (local) maximizer of $L_K$ in $\R^n$.

\noindent
The main result of this paper is the following strong version of the result of \cite{CCG}:
\noindent
\begin{theorem}
\label{th-1}
If a convex body $K$ in $\R^n$ is a local maximizer for $L_K$, then it has no positive
generalized Gauss
curvature at any point of its boundary. The same is true for  a centrally symmetric $K$
which is
a local maximizer for $L_K$ among centrally symmetric convex bodies.
\end {theorem}
\noindent
 An open problem is whether a maximizer for $L_K$ is
necessarily a polytope. Our result is a step in this direction,
because it shows that a maximizer has generalized Gauss curvature
equal to $0$ almost everywhere and never positive on its
boundary. \noindent To prove theorem \ref{th-1}, we shall suppose
that a convex body $K$ has a a positive generalized curvature at
some point $X_0$ of its boundary (see Definition \ref{def} below),
modify slightly $K$ in a neighborhood of $X_0$, from inside and
from outside to get a body $K'$ for which we shall estimate
$L_{K'}$. The paper is organized as follows. In section 2, after
presenting some notations, we study the effect of such modifications,
that are
 described in the general case in Lemmas \ref{lem-1}, \ref{lem-2} and \ref{lem-3},  and in the neighborhood
of some special points of the boundary of $K$ in Proposition \ref{prop-4} and Lemma \ref{lem-5}. Corollary \ref{cor-6} is a generalization of \cite{CCG}'s result, replacing positive curvature by strict convexity on an open subset of the boundary. To estimate carefully the asymptotic behavior of $L_{K'}$, we prove the geometric Lemma \ref{lem-7}  and we get in  Lemma \ref{lem-8} a special property of potential maximizers of $L_K$. Finally
section 3 is devoted to the proof of theorem \ref{th-1}, which needs some technical and very precise computations of volumes.

\noindent In connection to Theorem~\ref{th-1}, one should mention the paper \cite{RSW}, by Reisner, Sch\"utt and Werner,
where an analogous result is proved related to Mahler's conjecture. Namely: a minimizer $K$ of the volume-product can
not have a point of positive generalized Gauss curvature on its boundary (see also \cite{GM}).

\vskip 4mm\noindent
\section{Notations and preliminary results.}
Let $K$ be a convex body in $\R^n$.
It is not hard to show, and is well known, that for any convex body $K$, denoting by $g(K)$ the centroid of $K$, one has
$$M_K^{2n}=\frac{1}{n!}\int_{K-g(K)}\dots\int_{K-g(K)} \big(\det(X_1, \dots, X_n)\big)^2 dX_1\dots dX_n$$
$$=\frac{1}{n!}\int_{K}\dots\int_{K} \big(\det(X_1-g(K), \dots, X_n-g(K))\big)^2 dX_1\dots dX_n.$$
Let $X_0\in \partial K$. For $r>0$, denote $B(X_0,r)$  the Euclidean ball of center $X_0$ and radius $r$.

\begin{definition} \label{def}
 We say that $K$ has  positive  generalized (Gauss) curvature at $X_0$, if there exists  an inner normal
 $N$ of $K$ at $X_0$  and a positive definite quadratic form $q$ on $N^{\perp}=\{x\in \R^n; \langle x,N\rangle =0\}$
 such that
for every $\eps>0$, there exists $a>0$, such that
whenever $Y\in  N^{\perp}$ and $y\in \R$ satisfy
$$X_0+ Y+yN\in \partial K\cap B(X_0,a),$$ then
 $$(1-\eps)q(Y)\le y\le (1+\eps)q(Y).$$
 \end{definition}
 \noindent
 Of course, this normal $N$ and the quadratic form $q$ are then unique.  Observe that if $K$ is $C^2$ with positive
 curvature, then $K$ has positive generalized curvature at any point $X$ of its boundary, but that positive
 generalized curvature at some point $X_0$ does not imply any regularity  at any point of $\partial K$ other than $X_0$.
 We refer to \cite{SW} for more details on positive generalized curvature.
\vskip 1mm\noindent
The following two lemmas show the effect  of local slight modifications of an isotropic body $K$ on $\int_{K}\dots\int_{K} \big(\det(X_1, \dots, X_n)\big)^2 dX_1\dots dX_n$.
\vskip 4mm\noindent
\begin{lemma}
\label{lem-1}
Let $K$ be an isotropic convex body. Suppose that $C_m, m\ge 1$ is a sequence of Borel subsets of $\R^n$ such that
$C_m\cap {\rm int}(K)=\emptyset$,  $|C_m| >0$, $|C_m|\to 0$ and $K_m:= K\cup C_m$ is a convex body.
Then, when $m\to +\infty$,
$$\frac{1}{n!}\int_{K_m}\dots\int_{K_m} \big(\det(X_1, \dots, X_n)\big)^2 dX_1\dots dX_n
= M_K^{2n} +M_K^{2(n-1)}\int_{C_m} |X|^2 dX+ O(|K_m\setminus K|^2).$$
\end{lemma}

\begin{lemma}
\label{lem-2}
Let $K$ be an isotropic convex body. Suppose that  $D_m, m\ge 1$ is a sequence of Borel subsets of $\R^n$ such that $D_m\subset K$, $|D_m| >0$, $|D_m|\to 0$ and  $K'_m:= K\setminus D_m$ is a convex body.
Then, when $m\to +\infty$,
$$\frac{1}{n!}\int_{K'_m}\dots\int_{K'_m} \big(\det(X_1, \dots, X_n)\big)^2 dX_1\dots dX_n
= M_K^{2n}- M_{K}^{2(n-1)}\int_{D_m} |X|^2 dX+ O(|K\setminus K'_m|^2).$$
\end{lemma}
\vskip 4mm
\noindent
{\bf Proof of Lemma~{\ref{lem-1}}  and Lemma~{\ref{lem-2}}\,:}
\vskip 2mm
\noindent
One has
$$\frac{1}{n!} \int_{K_m}\dots \int_{K_m} \big(\det(X_1, \dots, X_n)\big)^2 dX_1\dots dX_n $$
$$=\frac{1}{n!} \Big(\int_{K}\dots \int_{K} \big(\det(X_1, \dots, X_n)\big)^2 dX_1\dots dX_n\ $$
$$+n\int_{C_m}\int_{K}\dots \int_{K} \big(\det(X_1, \dots, X_n)\big)^2 dX_1\dots dX_n+O(|K_m\setminus K|^2)\Big).$$
Now
$$\int_{C_m}\int_{K}\dots \int_{K} \big(\det(X_1, \dots, X_n)\big)^2 dX_1\dots dX_n$$
$$=\int_{C_m}\int_{K}\dots \int_{K}\big(\sum_{\sigma\in S_n}\sum_{\tau\in S_n}(-1)^{\varepsilon(\sigma)
\varepsilon(\tau)}   \prod_{i=1}^nX_{i\sigma(i)}X_{i\tau(i)}\big)dX_1\dots dX_n.$$
Since $K$ is isotropic one has
\begin{equation}
\label{proof-lem-1-eq-1}
\int_K X_{i\sigma(i)}X_{i\tau(i)}dX_i=0 \mbox{ if } \sigma(i)\not=\tau(i)\,.
\end{equation}
It follows that
$$\int_{C_m}\int_{K}\dots \int_{K} \big(\det(X_1, \dots, X_n)\big)^2 dX_1\dots dX_n$$
$$=\int_{C_m}\int_{K}\dots \int_{K}\big(\sum_{\sigma\in S_n}\sum_{\tau\in S_n}(-1)^{\varepsilon(\sigma)
\varepsilon(\tau)}   \prod_{i=1}^nX_{i\sigma(i)}X_{i\tau(i)}\big)dX_1\dots dX_n$$
$$=\sum_{\sigma\in S_n}\int_{C_m}\int_{K}\dots \int_{K} \prod_{i=1}^nX_{i\sigma(i)}^2 dX_1\dots dX_n$$
$$= (n-1)!\ M_K^{2(n-1)}\int_{C_m}\sum_{m=1}^n X_{1m}^2 dX_1$$
$$=  (n-1)!\ M_K^{2(n-1)}\int_{C_m} |X|^2 dX.$$
We can thus conclude. The proof of lemma 2 is analogous.\hskip 2mm $\square$
\vskip 2mm \noindent
In the next lemma, we investigate, under the hypotheses of lemmas \ref{lem-1} and \ref{lem-2} how $M_{K_m}$ differ from $M_K$.
\begin{lemma}
\label{lem-3}
 Under  the hypotheses of Lemma \ref{lem-1} or respectively of Lemma \ref{lem-2}, one has
$$M_{K_m}^{2n}=\frac{1}{n!}\int_{K_m}\dots\int_{K_m} \big(\det(X_1, \dots, X_n)\big)^2 dX_1\dots dX_n + O(|K_m\setminus K|^2)$$
or respectively,
$$M_{K'_m}^{2n}=\frac{1}{n!}\int_{K'_m}\dots\int_{K'_m} \big(\det(X_1, \dots, X_n)\big)^2 dX_1\dots dX_n + O(|K\setminus K'_m|^2). $$
\end{lemma}
\vskip 2mm \noindent
{\bf Proof:}
\quad We assume throughout the proof that $K$ is isotropic but, a posteriori, the equalities stated in the lemma
remain true under invertible linear transformations.

\noindent
Let  $g_m$ be the centroid of $K_m$. One has :
$$M_{K_m}^{2n}=\frac{1}{n!} \int_{K_m-g_m}\dots \int_{K_m-g_m} \big(\det(X_1, \dots, X_n)\big)^2 dX_1\dots dX_n.$$
Since the centroid of $K$ is at $0$. One has for every $u\in S^{n-1}$,
$$\langle g_m,u\rangle =\frac{1}{|K|+|C_m|} (\int_{K } \langle X,u\rangle dX+ \int_{C_m } \langle X,u\rangle dX)
=\frac{1}{|K|+|C_m|}\int_{C_m } \langle X,u\rangle dX\,,$$
and thus $|g_m|=O(|C_m|)$ (observe that the hypotheses imply that the $C_m$, $m\ge 1$, are uniformly bounded).

\noindent
We have
$$n!\ M_{K_m}^{2n}= \int_{K_m}\dots \int_{K_m} \big(\det(Y_1-g_m, \dots, Y_n-g_m)\big)^2 dY_1\dots dY_n$$
$$=\int_{K_m}\dots \int_{K_m} \big(\det(Y_1,\dots,Y_m)-\sum_{k=1}^n \det(Y_1,\dots,Y_{k-1}, g_m, Y_{k+1},\dots,Y_n)\big)^2 dY_1\dots dY_n$$
$$={\bf A}-{\bf B}+{\bf C}.$$
where
$${\bf A}:=\int_{K_m}\dots \int_{K_m} \big(\det(Y_1,\dots,Y_n)^2 dY_1\dots dY_n$$
$${\bf B}:=2\sum_{k=1}^n\int_{K_m}\dots \int_{K_m}\det(Y_1,\dots,Y_n)\det(Y_1,\dots,Y_{k-1}, g_m, Y_{k+1},\dots,Y_n) dY_1\dots dY_n$$
$${\bf C}:= \int_{K_m}\dots \int_{K_m}(\sum_{k=1}^n \det(Y_1,\dots,Y_{k-1}, g_m, Y_{k+1},\dots,Y_n)\big)^2 dY_1\dots dY_n$$
The term {\bf A} has been treated already :
$$\frac{{\bf A}}{n!}= M_K^{2n} +M_K^{2(n-1)}\int_{C_m} |X|^2dX +O(|C_m|^2).$$
Since $|g_m|=O(|C_m|)$, it is clear that
$${\bf C}=O(|C_m|^2).$$
For {\bf B} we write
$$\frac{{\bf B}}{2}={\bf D}+{\bf E}+O(|C_m|^2)$$
where
$${\bf D}:=\sum_{k=1}^n\int_{K}\dots \int_{K}\det(Y_1,\dots,Y_n) \det(Y_1,\dots,Y_{k-1}, g_m, Y_{k+1},\dots,Y_n) dY_1\dots dY_n$$
and
$${\bf E}:=\sum_{k=1}^n \int_{K}\dots\int_K\int_{C_m}\int_K \dots\int_{K} \det(Y_1,\dots,Y_n)\det(Y_1,\dots,Y_{k-1}, g_m, Y_{k+1},\dots,Y_n)dY_1\dots dY_n.$$
It is easily seen that ${\bf D}=0$, because of the isotropicity of $K$. Now, once again since
$g_m=O(|C_m|)$, one has ${\bf E}=O(|C_m|^2)$.
\vskip 0.5mm\noindent
The  corresponding result for $K'_m$ is proved in the same  way.\hskip 2mm$\square$

\begin{proposition}
\label{prop-4}
Under the assumptions of Lemma \ref{lem-1} or, respectively, Lemma~\ref{lem-2}
one has
\begin{equation}
\label{prop-4-eq-1}
L_{K_m}^{2n}=L_K^{2n}
\left[ 1+\frac{\int_{K_m\setminus K}|X|^2 dX  }{M_K^2}
-(n+2)\frac{|K_m\setminus K|}  {|K|}
 +O(|K_m\setminus K|^{2})\right]
 \end{equation}
 or, respectively,
 \begin{equation}
 \label{prop-4-eq-2}
 L_{K'_m}^{2n}=
 L_K^{2n}
\left[ 1-\frac{\int_{ K\setminus K'_m }|X|^2 dX  }{M_K^2}
+(n+2)\frac{|K\setminus K'_m|}  {|K|}
 +O(|K\setminus K'_m|^{2})\right]\,.
 \end{equation}
 \end{proposition}

 \vskip 1mm
\noindent
{\bf Proof:}
By Lemma~\ref{lem-1} and Lemma~\ref{lem-3} we have
$$L_{K_m}^{2n}= \frac{ M_K^{2n} +M_K^{2(n-1)}\int_{K_m\setminus K}
|X|^2 dX+ O(|K_m\setminus K|^2)}{|K|^{n+2}+ (n+2)|K|^{n+1}|K_m\setminus K|+O(|K_m\setminus K|^2)}\,.$$
From this (\ref{prop-4-eq-1}) follows. The equality (\ref{prop-4-eq-2}) is proved in a similar way.
\hskip 2mm
$\square$

\begin{lemma}
\label{lem-5}
Suppose that $K$ is an isotropic convex body and that,
 in addition to the conditions of Proposition~\ref{prop-4},
 there exists $X_0\in \partial K$ such that $X_0$ is in the closure of $C_m$ for all $m$ and\/  $\diam(C_m)\to 0$
 and also, $X_0$ is in the closure of $D_m$ for all $m$ and\/  $\diam(D_m)\to 0$. Then, if $K$ is a local maximizer or a local minimizer
 for $L_K$, we have
 \begin{equation}
 \label{lem-5-eq-1}
 |X_0|^2 |K|= (n+2) M_K^2\,.
 \end{equation}
\end{lemma}

\vskip1mm\noindent
{\bf Proof:}\ The conditions of the lemma imply that, when $m\to +\infty$, one has:
\begin{equation}
\label{lem-5-eq-2} \int_{K_m\setminus K}|X|^2 dX \sim
|X_0|^2|K_m\setminus K|
 \hbox{\ and\ }\int_{K\setminus K'_m}|X|^2 dX \sim |X_0|^2|K\setminus K'_m|\, .
\end{equation}
thus the result follows from Proposition~\ref{prop-4}.\hskip 2mm $\square$

\vskip 2mm\noindent
\noindent
{\bf Remarks}
\vskip 0.5mm\noindent
{\bf 1)} A common example of a point $X_0$ that satisfies the assumptions of Lemma~\ref{lem-5} is the following:
Let $X_0\in \partial K$. We say that {\em $\partial K$ is locally strictly convex at $X_0$\/} or that {\em $X_0$
is a point of local strict convexity of $\partial K$\/}, if there exists no non-degenerate line segment $I\subset
\partial K$ such that $X_0\in I$ (even as an end-point). The following claim is easy to prove:

\noindent
{\bf Claim.}\ {\em
Let $X_0$ be a point of local strict convexity of $\partial K$ and let $N\in S^{n-1}$ be an outer normal of $K$ at
$X_0$. Then the sets
$$C_m=\conv\Big(K\cup\big(X_0+\frac{1}{m}N\big)\Big)\setminus K$$
and
$$D_m=\{X\in K;\langle X, N\rangle\ge \langle X_0,N\rangle-\frac{1}{m}\big\}$$
satisfy the conditions of Lemma~5.\/}

\vskip 1mm
\noindent
{\bf 2)}
If $X_0 \in\partial K$ is a point of positive generalized curvature of $\partial K$ then it is a point of local
strict convexity and thus satisfies the conditions of Lemma~\ref{lem-5}.

\vskip 2mm
\noindent

As a corollary of Lemma~\ref{lem-5} and of \cite{CCG} (or of our Theorem~\ref{th-1}) we get the following strengthening
of a result of \cite{CCG}:

\begin{corollary}
\label{cor-6}
Suppose that there exists an open neighborhood $U$ in $\partial K$ which is strictly convex (that is, every point in $U$
is a point of local strict convexity). Then $K$ is not a local maximizer for $L_K$.
\end{corollary}
\noindent
{\bf Proof:}\quad We may assume that $K$ is isotropic.
By Lemma~\ref{lem-5} and the Claim following it, all the points in
$U$ have the same Euclidean norm. Thus $U$ is an open neighborhood
on a Euclidean sphere. The result of \cite{CCG} or
Theorem~\ref{th-1} now complete the proof.\hskip 2mm $\square$

\noindent
We shall later need the following geometric lemma.
\begin{lemma}
\label{lem-7}

Suppose that $K$ is a convex body containing $0$ in its interior and that $\partial K$ has
positive generalized curvature at some point $X_0$. Assume that the normal vector of $K$ at $X_0$ is not parallel to
the vector $X_0$. Then there exists $u\in S^{n-1}$ and $\alpha >0$ such that if
$K(\alpha,u)= \{X\in K;\langle X,u\rangle \ge\alpha\}$, then  $K(\alpha,u)$ is a cap of $K$  with non-empty interior
and $$\max_{X\in K(\alpha,u)} |X|<|X_0|.$$
\end{lemma}

\vskip 2mm
\noindent
{\bf Proof:}\quad After an affine change of variables in $\R^n$, transforming $0$ into $X_0$, we may suppose  that for $|Z|\le a$,
the boundary of $K$ is described by $z=g(Z)$ with $(Z,z)\in \R^n=\R^{n-1}\times \R$, and
$$(1-\eps)|Z|^2\le g(Z)\le (1+\eps)|Z|^2.$$
This affine change of variables transforms $B(0,|X_0|)$ into an ellipsoid ${\cal E}$ with $0\in \partial {\cal E}$,
 whose inner normal
$N$ at $0$ is not $e_n$. We may suppose that $N=\cos(\theta) e_1+\sin(\theta)e_n$ for some angle $\theta\in [0,\frac{\pi}{2}[$.
Also, since ${\cal E}$ has positive curvature at $0$, one can find some positive constants $b$ and $C$  such that
\begin{equation}
\label{paraboloid}
B(0,b)\cap{\cal P} \subset B(0,b)\cap{\cal E}
\end{equation}
where ${\cal P}$ is the paraboloid defined by
$${\cal P}=\{M:=xe_1 +Y+ze_n; \langle OM,N\rangle\ge  C( |OM|^2-\langle OM,N\rangle ^2)\}. $$
\vskip 0.5mm \noindent
Let $0<x_0<a$. The hyperplane $H$ tangent to the upper paraboloid
($z=(1+\eps)|Z|^2$) at $M_0=x_0e_1+(1+\eps)|x_0|^2 e_n$ has the equation
 $$z= (1+\eps)(2xx_0- x_0^2), $$
 where $M=xe_1 +Y+ze_n$ is a point in $\R^n$, with $Y\in \{e_1,e_n\}^{\perp}$.
 The zone ${\cal A}$ between the hyperplane $H$ and the lower paraboloid ($z=(1-\eps)|Z|^2$) is described by
 $${\cal A}=\{M: xe_1 +Y+ze_n;      (1-\eps) (x^2+ |Y|^2) \le              z\le (1+\eps)(2xx_0- x_0^2)\}$$
 Thus for $M\in {\cal A}$, one has
 $$x^2-2\frac{1+\eps}{1-\eps} x_0 x +\frac{1+\eps}{1-\eps} x_0^2\le 0$$
 which says that
 $$\Big(x-\frac{1+\eps}{1-\eps} x_0\Big)^2 \le \frac{1+\eps}{1-\eps}\Big(\frac{1+\eps}{1-\eps}-1\Big) x_0^2$$
 or
 $$     \Big(  \frac{1+\eps}{1-\eps} - \frac{\sqrt{2\eps(1+\eps)}} {1-\eps}  \Big)x_0   \le x \le  \Big(\frac{1+\eps}{1-\eps} +
 \frac{\sqrt{2\eps(1+\eps)}}   {1-\eps} \Big)x_0\ .$$
 It follows that for $\eps$ small enough one has for $M=xe_1+Y+ze_n\in {\cal A}$: $x < 2x_0$ and $x^2+ |Y|^2\le 3x_0^2$. Thus,
 for $x_0$ small enough,
 ${\cal A}\cap\{xe_1 +Y+ze_n; z\ge g(x,Y)\}$ is a cap of $K$.
 passing through $0$, with normal $N= \cos(\theta)e_1+\sin(\theta)e_n$.
 \vskip 1mm \noindent
 By (\ref{paraboloid}), it is sufficient to show that  for $x_0$ small enough, one has
 $${\cal A}\subset {\cal P}\cap B(0,b).$$
 First it is easy to choose $x_0$ small enough such that ${\cal A}\subset B(0,b)$
Observe then that
$${\cal P}=\{xe_1 +Y+ze_n;x\cos(\theta)+z\sin(\theta)\ge C\big(x^2+|Y|^2 -(x\cos(\theta)+z\sin(\theta))^2\big)\}$$
and that setting $x=x_0 u, Y=x_0 V$ and $z=x_0^2w$, one gets
$${\cal A}= \{x_0(u+V+x_0 w);  (1-\eps) (u^2+ |V|^2) \le w\le (1+\eps)(2u- 1)\}$$
Thus we need only to prove that if
 $(1-\eps) (u^2+ |V|^2) \le w\le (1+\eps)(2u- 1)$ then
 $$\frac{u\cos(\theta) +x_0 w\sin(\theta) }{ |V|^2+ (u\sin(\theta) +x_0 w\cos(\theta))^2}\ge Cx_0.$$
 which is clear when $x_0\to 0$ because $u\sim 1$ and $w$ is uniformly bounded.

\noindent
 Observe finally that if we have the singular case that the point of tangency $M=x_0e_1+ (1+\eps)|x_0|^2$
 of the upper paraboloid with the tangent hyperplane $H$ is on $\partial K$, then we get a cap of $K$
 by pushing $H$ a small distance into the upper paraboloid in the direction of its inner normal.
 \hskip 1mm$\square$


\vskip2mm\noindent

\begin{lemma}
\label{lem-8}
Under the assumptions of Lemma~\ref{lem-5}, if $K$ is a local maximizer for $L_K$ and $\partial K$ has
positive generalized curvature at $X_0$ then the outer normal $N(K,X_0)$ of $K$ at $X_0$ is parallel
to the vector $X_0$.
\end{lemma}

\vskip 1mm
\noindent
{\bf Proof:}
We assume that $L_K$ is maximal, $\partial K$ has positive generalized curvature at $X_0$ and the normal
vector of $K$ at $X_0$ is not parallel to $X_0$.

\noindent
Using Lemma~\ref{lem-7} we continue as follows: Let $u\in S^{n-1}$ and $\alpha>0$ be taken from Lemma~\ref{lem-7}.
Let $H=\{X;\, \langle X,u\rangle=\alpha\}$ and $H^+=\{X;\, \langle X,u\rangle\geq\alpha\}$.
Let $M=\max\{|X|;\, X\in H^+\cap K\}$. Then $M<|X_0|$.
Let  $d$ be the distance from $0$ to $H$, $h=h_K(u)-d$
and, for $m\geq 1$, let
$$D_m'=\{X\in K;\, h_K(u)-\frac{h}{m}\leq \langle X,u\rangle \leq h_K(u)\}\,.$$
Then the sequence $D_m'$ satisfies the conditions of Lemma~\ref{lem-2}. We have
$$\int_{D_m'}|X|^2\leq M^2|D_m'|\,.$$
Now, since $L_K$ is maximal, we have, combining the above with (\ref{prop-4-eq-2}), for $m$ big enough,
$$\frac{-M^2|D_m'|}{M_K^2}+(n+2)\frac{|D_m'|}{|K|}\leq O(|D_m'|^2)\,.$$
Combining the last inequality with (\ref{lem-5-eq-1}) we get, passing to the limit as $m\to \infty$,
$$|X_0|^2=\frac{(n+2)M_K^2}{|K|}\leq M^2 < |X_0|^2\,,$$
which is a contradiction.\hskip 2mm$\square$

\vskip 3mm
\noindent
\section
{\bf Proof of Theorem~\ref{th-1}.}

Assume that $K$ ia a local maximizer of $L_K$ and $X_0\in \partial K$ is a point of positive generalized  curvature
of $\partial K$. We may assume that $K$ is in isotropic position.

\noindent
By Lemma~\ref{lem-8}, we know that   $u=\frac{X_0}{|X_0|}$ is the external normal of $K$ at $X_0$.
We choose for $K_m$ and $K'_m$, $m\ge 1$,  the following sets:
$$K_m=\conv(X_0+\frac{u}{m}, K)$$
and
$$K'_m=\{X\in K;\ \langle X,u\rangle \le \langle X_0,u\rangle -\frac{1}{m}\}.$$

\vspace{2mm}
\noindent
By Remark~2) following Lemma~\ref{lem-5}, the sets $K_m\setminus K$ and $K\setminus K'_m$ satisfy the conditions
of Lemma~\ref{lem-5} and, of course, of Proposition~\ref{prop-4}.
In view of Lemma~\ref{lem-5}, it is essential to have an accurate estimation of
$$\int_{K_m\setminus K} |X|^2\,dX-|X_0|^2|K_m\setminus K|=\int_{K_m\setminus K}(|X|^2-|X_0|^2)\,dX$$
and
$$\int_{K\setminus K_m'} |X|^2\,dX-|X_0|^2|K\setminus K_m'|=\int_{K\setminus K_m'}(|X|^2-|X_0|^2)\,dX\,.$$

\noindent
For having such estimation it would be convenient to assume that the standard approximating ellipsoid of $K$
at $X_0$ is a Euclidean ball rather than just an ellipsoid.

\noindent
Let $u_1,\ldots,u_n$ be an orthonormal system in $\R^n$, with $u_n=\frac{X_0}{|X_0|}$ and such that
$u_1,\ldots,u_{n-1}$ are the directions of the principal radii of the quadratic form $q$ associated with $X_0$ (see Definition \ref{def}). Let
$T\in SL(n)$ be a volume preserving linear transformation of the form
$$T(\sum_{j=1}^n x_ju_j)=\sum_{j=1}^n \lambda_j x_ju_j\,; \quad \prod_{j=1}^n \lambda_j=1$$
(we write in short $T(X)=\Lambda X$ and $T^{-1}(X)=\Lambda^{-1}X$ assuming $X$ is written using the basis
$u_1,\ldots,u_n$). Choose $T$ so that the standard approximating ellipsoid of $\tilde K=T(K)$ at $T(X_0)$
is a Euclidean ball of radius $R$.

\noindent
Denoting $\tilde K_m=T(K_m)$ and $\tilde K_m'=T(K_m')$ we get
$$\int_{K_m\setminus K}(|X|^2-|X_0|^2)\,dX=\int_{\tilde K_m\setminus \tilde K}
(|\Lambda^{-1}Y|^2-|\Lambda^{-1}Y_0|^2)\,dY$$
and
$$\int_{K\setminus K_m'}(|X|^2-|X_0|^2)\,dX=\int_{\tilde K\setminus \tilde K_m'}
(|\Lambda^{-1}Y|^2-|\Lambda^{-1}Y_0|^2)\,dY\,.$$

\noindent
We shall use a temporary coordinate system that satisfies:
\begin{itemize}
\item[1)] $T(X_0)=0$
\item[2)] The outer normal vector of $\tilde K$ at $0$ is $-e_n$\hspace{1mm}($e_n$ is the $n$-th coordinate vector), thus
$\tilde K\subset \{X\in \R^n\,; \langle X,e_n\rangle \geq 0\}$
\end{itemize}
   \vskip 1mm\noindent
 We write $X=(Y,y)\in\R^n=\R^{n-1}\times \R$. Let $G=g(\tilde K)$ be the centroid of $\tilde K$. In our temporary
 coordinates $G=(0,b)$ with $b>0$ (in view of Lemma~\ref{lem-8}).
For $a>0$, small enough, define
$$C_a=\conv(\tilde K,(-a,0))\setminus \tilde K$$
$$D_a=\{(Y,y)\in \tilde K; y\le a\}.$$
By the above discussion, we have to estimate
for $\tilde K_m\setminus \tilde K=C_a$ and $\tilde K\setminus  \tilde K'_m=D_a$ \ ($a=\frac{1}{m}$),
the following quantities in terms of $a>0$, $a\to 0$:
$$\phi(a)=\int_{C_a} (|\Lambda^{-1}(X-G)|^2-|\Lambda^{-1}G|^2)dX$$
$$\psi(a)=\int_{D_a} (|\Lambda^{-1}(X-G)|^2-|\Lambda^{-1}G|^2)dX\,.$$
The equation of the boundary of the body, in a neighborhood of $0$ can be written as
$$y=\frac{|Y|^2}{2R}+ o(|Y|^2),$$

\noindent
With these notations
$$\phi(a)= \int_{(Y,y)\in C_a} \left(\sum_{j=1}^{n-1}\left(\frac{Y_j}{\lambda_j}\right)^2+
\left(\frac{y}{\lambda_n}\right)^2 -2\frac{yb}{\lambda_n^2}\right) dY dy,$$
$$\psi(a)= \int_{(Y,y)\in D_a} \left(\sum_{j=1}^{n-1}\left(\frac{Y_j}{\lambda_j}\right)^2+
\left(\frac{y}{\lambda_n}\right)^2 -2\frac{yb}{\lambda_n^2}\right) dY dy.$$

\noindent
We first estimate $\phi(a)$ and $\psi(a)$ under the hypothesis that in some neighborhood of $0$ the equation of the boundary
of $K$ is actually $$y=\frac{|Y|^2}{2R}.$$ Then we shall see that this approximation is actually good.

\vskip 2mm \noindent
{\bf 1)} We suppose that $y=\frac{|Y|^2}{2R}$. One has
$$D_a=\{ (Y,y)\in \R^n; |Y|\le \sqrt{2Ra},\ \frac{|Y|^2}{2R}\le y\le a\}\,.$$

Since $D_a$ is circular with respect to $Y$, we have
$$\int_{(Y,y)\in D_a}Y_j^2\,dY\,dy=\frac{1}{n-1}\int_{(Y,y)\in D_a}|Y|^2\,dY\,dy\,.$$
Substituting $\alpha_n=\frac{1}{n-1}\sum_{j=1}^{n-1}\lambda_j^{-1}$ we get
with a change of variable to polar coordinates in $\R^{n-1}$ and denoting  by $v_k$ the volume of the
Euclidean ball in $\R^k$,
$$\psi(a)=(n-1)v_{n-1} \int_{S_{n-2}}\int_0^{\sqrt{2Ra}}\big( \int_{\frac{r^2}{2R}}^a
(\alpha_n r^2+\lambda_n^{-1}(y^2 -2yb)) dy\big)r^{n-2}\,dr\,d\theta .$$
Setting  $r=\sqrt{2Ra}s$ and $y=az$ we get
$$\psi(a)=(n-1)v_{n-1}a(2Ra)^{\frac{n-1}{2}} \int_{S_{n-2}}\int_0^1\big(
\int_{s^2}^1 (2\alpha_n Ra s^2+\lambda_n^{-1}(a^2z^2 -2abz))\,dz\,s^{n-2}dsd\theta$$
$$= (n-1)v_{n-1}  (2R)^{\frac{n-1}{2}} a^{\frac{n+1}{2}}
\int_0^1 \big( 2\alpha_n Ras^{n}(1-s^2)+ \lambda_n^{-1}(\frac{1}{3}a^2 s^{n-2}(1-s^6)- ab(1-s^4)s^{n-2}) \big)ds$$
$$ = (n-1)v_{n-1}  (2R)^{\frac{n-1}{2}} a^{\frac{n+3}{2}}
\int_0^1 \big( 2\alpha_n Rs^{n}(1-s^2)+ \lambda_n^{-1}(\frac{a}{3} s^{n-2}(1-s^6)- b(1-s^4)s^{n-2}) \big)ds$$
$$= (n-1)v_{n-1}  (2R)^{\frac{n-1}{2}} a^{\frac{n+3}{2}}\Big((2\alpha_n R\big(\frac{1}{n+1}-\frac{1}{n+3}\big)-
\lambda_n^{-1}b\big(\frac{1}{n-1}-\frac{1}{n+3}\big) +O(a)\Big)$$
$$=4(n-1)v_{n-1}  (2R)^{\frac{n-1}{2}} a^{\frac{n+3}{2}}\Big( \frac{\alpha_n R}{(n+1)(n+3)}-
\frac{\lambda_n^{-1} b}{(n-1)(n+3)}\big)+O(a)\Big)$$
$$=\frac{4(n-1)v_{n-1}(2R)^{\frac{n-1}{2}}a^{\frac{n+3}{2}}}{(n+1)(n+3)}\cdot
\Big(\alpha_n R-\frac{(n+1)\lambda_n^{-1}b}{n-1}+O(a)
\Big).$$
We shall need also to compute $|D_a|$.
One has
$$|D_a|=  (n-1)v_{n-1}a(2Ra)^{\frac{n-1}{2}} \int_{S_{n-2}}\int_0^1 (1-s^2)s^{n-2}dsd\theta$$
$$=(n-1)v_{n-1}(2R)^{\frac{n-1}{2}}a^{\frac{n+1}{2}}(\frac{1}{n-1}- \frac{1}{n+1})
=\frac{2v_{n-1}}{n+1}(2R)^{\frac{n-1}{2}}a^{\frac{n+1}{2}}.$$

\vskip 2mm \noindent
{\bf 2)} We still suppose that the boundary of $\tilde K$ in a neighborhood of $0$ is given by $y=\frac{|Y|^2}{2R}$.
Then the tangent hyperplanes to $\tilde K$ through $(0,-a)$, indexed by  $\theta\in S^{n-2}$ - the direction of the projection
of their point of tangency with $\tilde K$, are given
by the equations
$$y=-a+\sqrt{\frac{2a}{R}}\langle \theta,Y\rangle\,.$$
It follows that
$$C_a=\Big\{ (Y,y)\in \R^n;\ |Y|\le \sqrt{2Ra},\ -a+\sqrt{\frac{2a}{R} } |Y|\le y\le\frac{|Y|^2}{2R}\Big\}$$
$$= \big\{(\sqrt{2Ra} Z,az)\in \R^n;\ |Z|\le 1,\ 2|Z|-1\le z\le |Z|^2\big\}.$$
Thus, using the same rotation invariance as in (1),
 $$\phi(a)= (n-1)v_{n-1}a(2Ra)^{\frac{n-1}{2}}\int_{S_{n-2}}\int_0^1
\big( \int_{2s-1}^{s^2} (2\alpha_n Ra s^2+\lambda_n^{-1}(a^2z^2 -2abz))s^{n-2}ds$$
$$=(n-1)v_{n-1}  (2R)^{\frac{n-1}{2}} a^{\frac{n+3}{2}}\Big(
\int_0^1 \big((2s^{n}(1-s)^2\alpha_n R -\lambda_n^{-1}b(s^4-(2s-1)^2) s^{n-2}\big) ds+O(a)\Big)$$
$$=(n-1)v_{n-1}  (2R)^{\frac{n-1}{2}} a^{\frac{n+3}{2}}\Big(\big( \frac{1}{n+1}-\frac{2}{n+2}+
\frac{1}{n+3}\big)2\alpha_n R-\big( \frac{1}{n+3}-\frac{4}{n+1}+\frac{4}{n}-\frac{1}{n-1}) \big)\lambda_n^{-1}b +O(a)\Big)$$
$$=(n-1)v_{n-1}  (2R)^{\frac{n-1}{2}} a^{\frac{n+3}{2}}\Big(\frac{4\alpha_n R}{(n+1)(n+2)(n+3)}
-\frac{4(n-3)\lambda_n^{-1}b}{n(n+1)(n-1)(n+3)} +O(a)\Big)$$
$$=\frac{4(n-1)v_{n-1}(2R)^{\frac{n-1}{2}}a^{\frac{n+3}{2}}}{(n+1)(n+3)}\cdot
\Big(\frac{\alpha_n R}{n+2}-\frac{(n-3)\lambda_n^{-1}b}{n(n-1)}+O(a)\Big).$$
Moreover
$$|C_a|= (n-1)v_{n-1}a(2Ra)^{\frac{n-1}{2}}\int_{S_{n-2}}(\int_0^1
(1-s)^2 s^{n-2}ds) d\theta= \frac{2v_{n-1}}{n(n+1)} (2R)^{\frac{n-1}{2}} a^{\frac{n+1}{2}}.$$
\vskip 2mm \noindent
{\bf 3)} But the hypothesis which has been done that in a neighborhood of $0$, the equation of the
boundary of $\tilde K$ is  $y=\frac{|Y|^2}{2R}$ has to be replaced with the following one: For every $\eps>0$, there exists $c\ge 0$ such that
$$(1-\eps)\frac{|Y|^2}{2R}\le  y\le (1+\eps)\frac{|Y|^2}{2R}\hbox{ whenever } |Y|^2+y^2\le c.$$
One has to see that in terms of $a$, the estimates of {\bf 2)} and {\bf 3)} still hold.
We shall treat first $\psi(a)$ and then $\phi(a)$.
\vskip 1mm \noindent
One has $$D_a\subset
\{ (Y,y)\in \R^n; |Y|\le \sqrt{2R_+(a)a},\ \frac{|Y|^2}{2R_+(a)}\le y\le a\}\,,$$
$$ \{ (Y,y)\in \R^n; |Y|\le \sqrt{2R_-(a)a},\ \frac{|Y|^2}{2R_-(a)}\le y\le a\}\subset D_a$$
 and
$$C_a\subset\{(Y,y)\in\R^n;\,|Y|\leq \sqrt{2R_+(a)a},\,-a+\sqrt{\frac{2a}{R_+(a)}}\leq y
\leq\frac{|Y|^2}{2R_-(a)}\}\,,$$
$$\{(Y,y)\in\R^n;\,|Y|\leq \sqrt{2R_-(a)a},\,-a+\sqrt{\frac{2a}{R_-(a)}}\leq y
\leq\frac{|Y|^2}{2R_+(a)}\}\subset C_a\,,$$
with $R_+(a)=R+\eps_+(a)$ and $R_-(a)= R-\eps_-(a)$, where $\eps_+(a)$ and $\eps_-(a)$ are nonnegative functions tending to $0$ when $a\to 0$.
Then everything works with upper and lower bounds for the negative and the positive terms on $D_a$ and $C_a$,
observing also that that $|D_a|^2$ and $|C_a|^2$ are of the order of $a^{n+1}$ which is negligible  with respect to $a^{\frac{n+3}{2}}$,
so that we can apply Proposition \ref{prop-4}.

\vskip 2mm
\noindent
{\bf Remark.}\quad The importance of Lemma~\ref{lem-8} comes in step {\bf 3)} above. Here, if the normal vector of
$\tilde K$ at $0$ were not parallel to the $y$-axis, we would get an extra error term of order that could be estimated only by
$a^{\frac{n+2}{2}}o(a)$. For our proof of Theorem~\ref{th-1} to work we would need an estimate of order
$a^{\frac{n+3}{2}}o(a)$ for this term.

\vskip 2mm
\noindent
To conclude, using Proposition~\ref{prop-4} and Lemma~\ref{lem-5} (including (\ref{lem-5-eq-2}) in its proof) and
replacing $K_m\setminus K$ by $T^{-1}(C_a)$ and
  $K\setminus K'_m$ by $T^{-1}(D_a)$, the above computations show that
for some functions $c(n,R)$ and $d(n,R)$ depending only of $n$ and $R$,
  $$L_{K_m}^{2n}=L_K^{2n}\Big(1+ c(n,R) a^{ \frac{n+3}{2}}\Big(\alpha_n R
-\frac{(n+2)(n-3)}{n(n-1)}\lambda_n^{-1}b +O(a)\Big)\Big)$$
and
$$L_{K'_m}^{2n}= L_K^{2n}\Big(1-d(n,R) a^{ \frac{n+3}{2}}
\Big( \alpha_n R- \frac{n+1}{n-1}\lambda_n^{-1}b+O(a)\Big)\Big).$$
Thus one has both
$$\alpha_n\lambda_n R
\le \frac{(n+2)(n-3)}{n(n-1)}b \hbox{\hskip 3mm and \hskip 2mm } \alpha_n\lambda_n R\ge \frac{n+1}{n-1} b,$$
So that $$\frac{(n+2)(n-3)}{n(n-1)}\ge \frac{n+1}{n-1}$$
which gives a contradiction.
 \vskip 2mm \noindent
 Note that in the case that $K$ is centrally symmetric, a similar argument, using $C_m$ and $-C_m$ together
 and $D_m$ and $-D_m$ together will work in the same way, keeping $K_m$ and $K'_m$ centrally symmetric.
 This observation takes care of the centrally symmetric part of Theorem~\ref{th-1}. There the use of
 lemma \ref{lem-3} is not needed, due to symmetry. \hskip 2mm$\square$
\vskip 3mm \noindent
{\bf Acknowledgements:}\quad
Mathieu Meyer and Shlomo Reisner are grateful respectively to the University of Haifa and the Technion and to the LAMA at Univ.\ Paris-Est Marne-la-Vall\'ee for their hospitality and support  during part
of the work on this paper. Both thank Yehoram Gordon for helpful discussions with him.

\vskip 3mm \noindent
Mathieu Meyer \newline
Universit\'e Paris-Est\newline
Laboratoire d'Analyse et de Math\'ematiques Appliqu\'ees (UMR 8050) UPEMLV\newline
F-77454 Marne-la-Vall\'ee Cedex 2, France\newline
{\tt mathieu.meyer@u-pem.fr}

\vskip 2mm \noindent
Shlomo Reisner \newline
Department of Mathematics \newline
University of Haifa \newline
Haifa, 31905, Israel \newline
{\tt reisner@math.haifa.ac.il}

\end{document}